\documentclass{article}

\usepackage{amssymb}
\usepackage{amsfonts}
\usepackage{amsmath}
\usepackage{amsthm}
\usepackage[numbers]{natbib}

\usepackage{color,mathdots,setspace,framed,graphicx,natbib,latexsym,nicefrac,longtable}

\title{A Topological Tennenbaum Theorem}
\author{Elliot Glazer}
\date{November 20, 2023}

\begin{document}

\maketitle

This paper provides negative answers to all parts of the following question that was posed by Enayat, Hamkins, and Wcis\l{}o as Question 3.6 of \cite{Enayat-Hamkins-Wcislo}.\footnote{The origin of the question goes back to an unpublished note of Enayat that was circulated among the participants of a meeting in 2009 at the Mittag-Leffler Institute.}  

\medskip

\noindent \textbf{Question.~}\emph{Does the Baire space support a topological model of $\mathrm{PA}$?
Does $\mathbb{R}^{\omega}$ support a topological model of $\mathrm{PA}$? Does any uncountable Polish space support a topological model of $\mathsf{PA}$}?

\medskip

In the above question, $\mathrm{PA}$ is (first order) Peano Arithmetic, and a topological model of arithmetic is one in which the addition and multiplication operations are continuous (relative to a given topological space). As explained in Remark 3.7  of \cite{Enayat-Hamkins-Wcislo}, the problem was motivated by a version of Borel completeness theorem (for first order logic) established by Malitz, Mycielski, and Reinhardt  \cite{Malitz-Mycielski-Reinhardt}.

\medskip

 The aforementioned result of Malitz, Mycielski, and Reinhardt states that every countable theory $T$ with an infinite model has a model with universe the Baire space such that all definable relations are $\mathbf{\Delta}^0_2$ with respect to the Baire space topology. This is a topological analogue of the arithmetized completeness theorem that states that $T$ has a model with universe $\omega$ such that every definable relation is $\Delta^{0,T}_2$ in the arithmetical hierarchy relative to the oracle $T$. The latter fact is demonstrated to be sharp by Tennenbaum's Theorem, that in a
 nonstandard model of $\mathsf{PA}$ (Peano Arithmetic) $(\omega, \oplus, \otimes),$ neither
 $\oplus$ nor $\otimes$ is computable. Theorem 1 analogously demonstrates
 MMR's Theorem to be sharp, and Theorem 2 generalizes this result to a topological version of Tennenbaum's Theorem, namely that in a Borel model of arithmetic on an uncountable Polish space, neither addition nor multiplication is continuous.

 \medskip
 
We work in $\mathrm{Z}_2$ (full second-order arithmetic) for now and discuss optimal metatheory towards the end of the paper.

\medskip

\noindent\textbf{Theorem 1.}
\emph{There is no model $\mathcal{M}=(M, +, \cdot) \models \mathrm{PA}$ such that $M$ is a perfect Polish space and $+$ and $\cdot$ are continuous}.

\medskip

\noindent \textbf{Proof:} Suppose towards contradiction there is such $\mathcal{M}.$ Let $d$ be a complete metric on $M$ and $\langle a_i: n<\omega \rangle$ enumerate a dense subset of $M.$ Notice that the order $<$ on $M$ is Borel, since $<$ and $\ge$ have analytic definitions.

\medskip

Let $x_n > 5$ be such that $d(0, x_n) < \frac{1}{n}.$ 
Set $y_0 = x_0$ and recursively define $y_{n+1} = x_k y_n,$
where $k$ is least such that for all $S \subset (n+1),$

$$d \left (\sum_S y_i, x_k y_n + \sum_S y_i \right ) < 5^{-n}
\min_{T \neq U \subset (n+1)} d \left (\sum_T y_i, \sum_U y_i \right).$$

\noindent For $r \in 2^{\omega}$, define $$z_r = \lim_N \sum_{r(n)=1, n<N} y_n.$$

It is routine to verify that for each $r,$ the partial sums form a Cauchy sequence, and that $r \mapsto z_r$ is injective. Furthermore, an inductive argument using continuity of addition shows that for disjoint $A, B \subset \omega$ with $B$ finite,
$$z_{\chi_A} + z_{\chi_B} = z_{\chi_{A \cup B}}.$$

Choose $r \in 2^{\omega}$ uniformly at random. Let $r' \in 2^{\omega}$ denote the complement, i.e. the map $n \mapsto 1-r(n).$ From symmetry, we have the probability $$P(z_{r'} < z_{r}) = \frac{1}{2},$$ noting that equality is impossible. Let $$d_r = |z_{r} - z_{r'}|.$$

We next show that almost surely, $z_r$ and $z_{r'}$ are far apart:

\medskip

\noindent\textbf{Claim.} For any fixed $n$ and $k$,
$$P(z_{r'} < z_{r} \wedge d_r < y_n) < 2^{-k}.$$

\medskip

\noindent \textbf{Proof of Claim:}
Let $E= \{r: z_{r'} < z_{r} \wedge d_r < y_n\}.$
It suffices to show that $r \mapsto r \restriction (\omega \setminus \{n+1, ..., n+k\})$ is injective on $E.$ Suppose not. Fix $r_1 \neq r_2 \in E$ which agree outside of $\{n+1,..., n+k\}$ with $d_{r_1} \ge d_{r_2}.$ Using the fact that $y_{i+1} > 5y_i$ for all $i$ (and identifying $r_i$ with $r_i^{-1}(\{1\})$):

\begin{align*}
d_{r_1} - d_{r_2} &= (z_{r_1} - z_{r_1'}) - (z_{r_2} - z_{r_2'})
\\&= \left (z_{r_1 \cap r_2} + \sum_{i \in r_1 \setminus r_2} y_i\right )
- \left (z_{r_1' \cup r_2'} - \sum_{i \in r_1 \setminus r_2} y_i \right )
\\&- \left (z_{r_1 \cap r_2} + \sum_{i \in r_2 \setminus r_1} y_i \right )
+ \left (z_{r_1' \cup r_2'} - \sum_{i \in r_2 \setminus r_1} y_i \right )
\\&= 2 \left ( \sum_{i \in r_1 \setminus r_2} y_i - \sum_{i \in r_2 \setminus r_1} y_i \right )> y_{\max_{r_1 \setminus r_2}} \ge y_n,
\end{align*}

\noindent so $d_{r_1} > y_n,$ contradiction. The Claim is proven. $\square$

\medskip

Applying symmetry and $\sigma$-additivity, the Claim implies
$$P(\forall n (z_{r'} + y_n < z_{r}))=P(z_{r'} < z_{r} \wedge  \forall n ( d_r > y_n))=\frac{1}{2}.$$

\noindent But this is a tail event, so this violates Kolmogorov's zero-one law, contradiction. This concludes the proof of the theorem. $\square$

\medskip

For our topological analogue of Tennenbaum's Theorem, we use the weakest hypotheses with which we know how to execute the above proof strategy. In particular, we will consider models of arithmetic on an arbitrary uncountable Polish space. The main difficulty this presents is that 0 may be an isolated point, so we need an alternate approach to finding a sequence $\langle y_n \rangle$ which increases exponentially in terms of arithmetic but whose partial sums over any subset are rapidly convergent under the metric. The key step is constructing terms $c$ and $h_{k,n}$ for which we can prove Claim 1: Subclaim 1
below.

This work-around is sufficiently flexible to allow us to achieve a contradiction from just one of $+$ and $\cdot$ being continuous, and the other being Borel. We set up Claim 1 (discontinuity of addition) and Claim 2 (discontinuity of multiplication) to have roughly parallel proofs.

\medskip

\noindent\textbf{Theorem 2.}
\emph{Let $\mathcal{M}=(M, +, \cdot) \models \mathsf{Q}+\mathrm{I}\Delta_0 + \mathsf{Exp}$ be such that $M$ is an uncountable Polish space and $+$ and $\cdot$ are Borel functions. Then neither $+$ nor $\cdot$ is continuous.}

\medskip

\noindent Note that we are using the formulation of $\mathsf{Exp}$ as a first-order axiom about addition and multiplication (``for all $x$ and $y,$ there is $z$ and $s$ such that $s$ is a G\"{o}del code for a sequence witnessing
$x^y=z"$).

We will only use $\mathsf{Exp}$ to justify existence of the unary function $\exp(x)=2^x,$ which in turn will only be used in verifying discontinuity of multiplication. In particular, our argument will immediately give the following:

\medskip

\noindent\textbf{Corollary 1.}
\emph{There is no model $\mathcal{M}=(M, +, \cdot)$ of $\mathsf{Q}+\mathrm{IOpen}$ such that $M$ is an uncountable Polish space, $+$ is continuous, and $\cdot$ is Borel}.\footnote{IOpen is the induction schema for quantifier-free formulas. It is weaker than $\mathrm{I}\Delta_0.$}

\medskip

\noindent \textbf{Proof:} We begin with some observations about $\mathcal{M}$ that depend only on $+$ and $\cdot$ being Borel. Let $\langle a_i \rangle \subset M$ enumerate a dense subset and let $d: M^2 \rightarrow \mathbb{R}$ be a complete metric. The order $<$ on $M$ is Borel, since $<$ and $\ge$ have analytic definitions. The graph of the binary function
$(x, y) \mapsto x \bmod y$ has an analytic definition, and is therefore Borel. We can define the graph of $\exp$ as
\begin{align*}
\exp=\{(x,y): &\exists (z_1, z_2)((z_1 \bmod (z_2+1) = 1)\wedge 
(z_1 \bmod ((x+1)z_2+1) = y)
        \\&\wedge \forall i < x (2(z_1 \bmod ((i+1)z_2+1))=
        (z_1 \bmod ((i+2)z_2+1)))\}.
\end{align*}

This is a $\mathbf{\Sigma^1_2}$ definition, which immediately implies the graph is $\mathbf{\Delta^1_2}.$ 

\medskip

\noindent \textbf{Lemma.} \emph{For every uncountable Borel
$X \subset M,$ there is an increasing embedding
$f: 2^{\omega} \rightarrow X.$}

\medskip

\noindent \textbf{Proof of Lemma:} By the perfect set property, there
is some embedding $g: 2^{\omega} \rightarrow X.$ Define a
Borel coloring $b: [2^{\omega}]^2 \rightarrow 2$ by
$b(\{s, t\}) = 1$ for $s < t$ if $g(s) < g(t)$. Galvin's theorem\footnote{Galvin's Borel version of Ramsey's theorem for pairs was established in his unpublished work in the 1960s; it was later generalized to higher exponents by Blass \cite{Blass-perfect_Ramsey}.} implies that $b$ is homogeneous on some perfect $P \subset 2^{\omega}$, from which we can extract an increasing embedding $f: 2^{\omega} \rightarrow X.$ The Lemma is proven. $\square$

\medskip

Consider the cut $$N = \bigcup_{|[0, x]| \le \aleph_0} [0, x].$$ Note that for any increasing embedding $f: 2^{\omega} \rightarrow M,$ $f^{-1}(N) \subset \{0\}.$ In particular\footnote{This can also be seen by a Fubini argument}, $N \subsetneq M.$


Fix a point $c \in M \setminus N$. Notice that for standard $n$, $c^{n^{-1}}$ (floor symbols suppressed) is also not in $N$.
Let
$$ X_n = \{c^{2 - 2^{-n}} i : i < c^{2^{-n-2}}  \}  $$

\noindent Each $X_n$ is the image of an uncountable Borel set under the injective Borel map $x \mapsto c^{2 - 2^{-n}} x,$ so $X_n$ is an uncountable Borel set. Applying the Lemma and $\Pi^1_1-\mathrm{AC}_0,$ there is a sequence $\langle f_n \rangle$ of increasing embeddings
$f_n: 2^{\omega} \rightarrow X_n.$

We now verify the discontinuity of the operations.

\medskip

\noindent \textbf{Claim 1:} The function $+$ is not continuous.

\medskip

\noindent \textbf{Proof of Claim:}
Suppose towards contradiction $+$ is continuous. 

\noindent Let $x_n=f_n(0)$ and
$$h_{k, n} = f_n(\chi_{\{k\}})-x_n \in [c^{2 - 2^{-n}}, \text{ } c^{2 - 2^{-n} + 2^{-n-2}}).$$

\medskip

\noindent \textbf{Subclaim 1:} If $m \ge c^2,$ then for any $n,$  we have $m + h_{k,n} \rightarrow_k m.$

\medskip

\noindent \textbf{Proof of Subclaim:}
Since $m \ge c^2 \ge x_n,$
\begin{align*}
m &= (m - x_n) + x_n = (m - x_n) + \lim_k  (f_n(\chi_{\{k\}})) = (m-x_n)+ \lim_k (x_n + h_{k,n})
\\&=\lim_k (m+h_{k,n}).
\end{align*}
\noindent The Subclaim is proven. $\square$

\medskip

We construct a sequence $\langle k_n \rangle$ by recursion. Let $k_0=0.$ Suppose we have constructed $\langle k_i: i<n \rangle$ for some $n>0.$ For $S \subset n,$ let
$\hat{h}_S = c^2 + \sum_{i \in S} h_{k_i, i}.$
Then let $k_n$ be least such that
$$\max_{S \subset n} d(\hat{h}_S, \hat{h}_S + h_{k_n, n} ) < 5^{-n} \min_{T \neq U \subset n} d(\hat{h}_T, \hat{h}_U). $$

\medskip

Let $y_n = h_{k_n, n} \ge c.$ Since $h_{k, n} \in [c^{2 - 2^{-n}}$, $c^{2 - 2^{-n} + 2^{-n-2}})$, we have for every $n$, $y_{n+1} > 5y_n$.
We now proceed as in the proof of Theorem 1.

For $r \in 2^{\omega}$, define $$z_r = \lim_N \left (c^2+\sum_{r(n)=1, n<N} y_n \right).$$

It is routine to verify that for each $r,$ the terms of this limit form a Cauchy sequence, and that $r \mapsto z_r$ is injective.

Choose $r \in 2^{\omega}$ uniformly at random. Let $r' \in 2^{\omega}$ denote the complement. From symmetry, we have the probability $$P(z_{r'} < z_{r}) = \frac{1}{2},$$ noting that equality is impossible. Let $$d_r = |z_{r} - z_{r'}|.$$

\medskip

\noindent\textbf{Subclaim 2.} For any fixed $n$ and $k$,
$$P(z_{r'} < z_{r} \wedge d_r < y_n) < 2^{-k}.$$

\medskip

\noindent \textbf{Proof of Subclaim:}
Let $E= \{r: z_{r'} < z_{r} \wedge d_r < y_n\}.$
It suffices to show that $r \mapsto r \restriction (\omega \setminus \{n+1, ..., n+k\})$ is injective on $E.$ Suppose not. Fix $r_1 \neq r_2 \in E$ which agree outside of $\{n+1,..., n+k\}$ with $d_{r_1} \ge d_{r_2}.$ Using the fact that $y_{i+1} > 5y_i$ for all $i$:

$$d_{r_1} - d_{r_2} = 2 \left ( \sum_{i \in r_1 \setminus r_2} y_i - \sum_{i \in r_2 \setminus r_1} y_i \right ) > y_{\max_{r_1 \setminus r_2}} \ge y_n,$$

\noindent so $d_{r_1} > y_n,$ contradiction. The Subclaim is proven. $\square$

\medskip

Applying symmetry and $\sigma$-additivity, the Claim implies
$$P(\forall n (z_{r'} + y_n < z_{r}))=P(z_{r'} < z_{r} \wedge  \forall n ( d_r > y_n))=\frac{1}{2}.$$

\noindent But this is a tail event, so this violates the zero-one law, contradiction. Thus, $+$ is not continuous, proving the Claim. $\square$

\medskip

\noindent \textbf{Claim 2:} The function $\cdot$ is not continuous.

\medskip

\noindent \textbf{Proof of Claim:}
Suppose towards contradiction $\cdot$ is continuous. Define a partial map $g: \omega^3 \rightharpoonup 2^{\omega}$
by letting $g(i, k,n) = s$ if $$(exp \circ f_n)^{-1}(B_{1/k}(a_i))=\{s\}.$$

Let $F: 2^{\omega} \rightarrow 2^{\omega}$ be an increasing embedding
such that $$\mathrm{rng}(F) \cap \mathrm{rng}(g) = \emptyset.$$

\noindent Let $F_n = f_n \circ F.$ Notice that for all $n,$ $(\exp \circ F_n): 2^{\omega} \rightarrow M$ is an injection such that the range $R_n:=\mathrm{rng}(\exp \circ F_n)$ has no isolated points.

Let $x_n= F_n(0).$ Use $\Sigma^1_2-\mathrm{AC}_0$ to choose
$s_{k,n} \in 2^{\omega} \setminus \{0\}$ such that 
$$d(\exp(x_n), \exp(F_n(s_{k,n})))< 2^{-k}.$$ Let $$h_{k,n}=F_n(s_{k,n})-F_n(0)\in [c^{2 - 2^{-n}}, c^{2 - 2^{-n} + 2^{-n-2}}).$$

\medskip

\noindent \textbf{Subclaim 1:} If $m \ge c^2,$ then for any $n,$
we have $\exp(m + h_{k,n}) \rightarrow \exp(m).$

\medskip

\noindent \textbf{Proof of Subclaim:}
Since $m \ge c^2 \ge x_n,$
\begin{align*}
\exp(m) &= \exp(m - x_n) \exp(x_n) = \exp(m - x_n) \lim_k \exp(F_n(s_{k,n}))
\\&= \exp(m-x_n) \lim_k \exp (x_n + h_{k,n})
\\&= \lim_k (\exp(m-x_n) \exp(x_n+h_{k,n}))
\\&=\lim_k \exp(m+h_{k,n}).
\end{align*}
\noindent The Subclaim is proven. $\square$

\medskip

\medskip

\noindent The proof will now proceed parallel to that of Claim 1.

We construct a sequence $\langle k_n \rangle$ by recursion. Let $k_0=0.$ Suppose we have constructed $\langle k_i: i<n \rangle$ for some $n>0.$ For $S \subset n,$ let
$$\hat{h}_S = \exp(c^2) \prod_{i \in S} \exp(h_{k_i, i}).$$
Then let $k_n$ be least such that
$$\max_{S \subset n} d(\hat{h}_S, \hat{h}_S \exp(h_{k_n, n}) ) < 5^{-n} \min_{T \neq U \subset n} d(\hat{h}_T, \hat{h}_U). $$

\medskip

Let $y_n = \exp(h_{k_n, n}) \ge \exp(c).$ Since $h_{k, n} \in [c^{2 - 2^{-n}}$, $c^{2 - 2^{-n} + 2^{-n-2}})$, we have for every $n$, $y_{n+1} > y_n^5.$
For $r \in 2^{\omega}$, define $$z_r = \lim_N \left (\exp(c^2) \prod_{r(n)=1, n<N} y_n \right ).$$

It is routine to verify that for each $r,$ the terms of this limit form a Cauchy sequence, and that $r \mapsto z_r$ is injective.

Choose $r \in 2^{\omega}$ uniformly at random. Let $r' \in 2^{\omega}$ denote the complement. From symmetry, we have the probability $$P(z_{r'} < z_{r}) = \frac{1}{2},$$
noting that equality is impossible. Let $$d_r = \left \lfloor z_{r} z_{r'}^{-1} \right \rfloor.$$

\medskip

\noindent\textbf{Subclaim 2.} For any fixed $n$ and $k$,
$$P(z_{r'} < z_{r} \wedge d_r < y_n) < 2^{-k}.$$

\medskip

\noindent \textbf{Proof of Subclaim:}
Let $E= \{r: z_{r'} < z_{r} \wedge d_r < y_n\}.$
It suffices to show that $r \mapsto r \restriction (\omega \setminus \{n+1, ..., n+k\})$ is injective on $E.$ Suppose not. Fix $r_1 \neq r_2 \in E$ which agree outside of $\{n+1,..., n+k\}$ with $d_{r_1} \ge d_{r_2}.$ Using the fact that $y_{i+1} > y_i^5$ for all $i$:

$$d_{r_1} d_{r_2}^{-1} \ge \frac{z_{r_1}z_{r'_2}}{2z_{r'_1}z_{r_2}}
= \frac{1}{2}\left ( \prod_{i \in r_1 \setminus r_2} y_i \prod_{i \in r_2 \setminus r_1} y_i^{-1} \right )^2  > y_{\max_{r_1 \setminus r_2}} \ge y_n,
$$

\noindent so $d_{r_1} > y_n,$ contradiction. The Subclaim is proven. $\square$

\medskip

Applying symmetry and $\sigma$-additivity, Subclaim 2 implies
$$P(\forall n (y_n z_{r'}   < z_{r}))=P(z_{r'} < z_{r} \wedge  \forall n ( d_r > y_n))=\frac{1}{2}.$$

\noindent But this is a tail event, so this violates the zero-one law, contradiction. Thus, $\cdot$ is not continuous, proving the Claim and the Theorem as well. $\square$

\medskip






\noindent \textbf{Optimization.} Theorem 1 and Corollary 1 can each be proven in $\mathsf{ATR}_0,$ the weakest theory that can reasonably speak about Borel sets. This theory is sufficient to develop the theory of measure (\cite{yu}) and
category (\cite{baire}) of Borel sets, as well as Suslin's theorem that $\mathbf{\Delta}^1_1$ sets are Borel (\cite{Simpson-SOSA}, \S V.3).

\medskip

The use of $\Pi^1_1-\mathrm{AC}_0$ in choosing the sequence of increasing embeddings $\langle f_n \rangle$ can be circumvented by direct construction.
From Borel codes for $\cdot$ and $<,$ and from the sequence $\langle c^{2^{-n-2}} \rangle,$ construct codes $b_n$ for $X_n.$ Use a pseudohierarchy (as in \cite{Simpson-SOSA}, \S V.4) to uniformly construct embeddings $g_n: 2^{\omega} \rightarrow X_n.$ Define Borel colorings $p_n: [2^{\omega}]^2 \rightarrow 2$ by $p_n(\{s, t\})=1$ for $s < t$ if $g_n(s)<g_n(t).$

Since $\bigsqcup_n p_n^{-1}(1)$ has the property of Baire, there are regular open $U_n \subset [2^{\omega}]^2$ and $\langle C_i \rangle$ closed nowhere dense sets such that for all $n,$
$$U_n \triangle p_n^{-1}(1) \subset \bigcup_i C_i.$$

Follow \cite{mycielski} to construct a perfect tree $P \subset 2^{<\omega}$ with
$$[P]^2 \cap \bigcup_i C_i =\emptyset.$$

\noindent The colorings $p_n \restriction [P]^2$ are open. Follow \cite{Blass-perfect_Ramsey} to construct perfect $P_n \subset P$ on which $p_n$ is constant, from which we easily extract increasing embeddings
$f_n: 2^{\omega} \rightarrow X.$

\medskip

\noindent \textbf{Further work.}
A large amount of comprehension is used in constructing the partial map $g$ in Claim 2, as well as to justify the use of $\Sigma^1_2-\mathrm{AC}_0$ in choosing $\langle s_{k,n}\rangle.$ We thus ask:

\medskip

\noindent\textbf{Question 1.} Is Theorem 2 provable in $\mathsf{ATR}_0$?
\medskip

Little is known about Polish models of weaker arithmetics. We end with the natural next point of investigation:

\medskip

\noindent\textbf{Question 2.} Does some uncountable Polish space support a model of Presburger arithmetic with continuous addition?

\medskip

We speculate that an affirmative answer to Question 1 would lead to a negative answer to Question 2, by analogy between circumventing use of the $\mathbf{\Delta}^1_2$ function $\exp$ in showing discontinuity of $\cdot$ and in circumventing use of multiplication in showing discontinuity of $+.$

\medskip

\noindent \textbf{Acknowledgment} I am grateful to Ali Enayat for showing me this problem and for many helpful discussions and suggestions on earlier drafts of this paper.

\bibliographystyle{plain}
\bibliography{Polish.bib}
\end{document}